\documentclass[10pt,reqno]{amsart}
\usepackage{amsmath,amsopn,amssymb,amsthm}
\usepackage[cp1251]{inputenc}
\usepackage[english]{babel}
\usepackage[pagebackref,breaklinks=true,colorlinks=true,linkcolor=blue,citecolor=blue,urlcolor=blue]{hyperref}
\usepackage{graphicx}%
\usepackage{color}

\newtheorem{theorem}{Theorem}[section]

\newtheorem{definition}[theorem]{Definition}

\newtheorem{lemma}[theorem]{Lemma}

\newtheorem{question}[theorem]{Question}

\newtheorem{remark}[theorem]{Remark}

\begin{document}

\title[Geodesic orbit and weak symmetric spray manifolds]{Geodesic orbit and weakly symmetric spray manifolds}
\author{Xiyun Xu}
\address[Xiyun Xu] {School of Mathematical Sciences,
Capital Normal University,
Beijing 100048,
P.R. China}
\email{2210502058@cnu.edu.cn}
\author{Ming Xu}
\address[Ming Xu] {Corresponding author, School of Mathematical Sciences,
Capital Normal University,
Beijing 100048,
P.R. China}
\email{mgmgmgxu@163.com}

\begin{abstract}
In this paper, we introduce the geodesic orbit and weakly symmetric properties in homogeneous spray geometry. When a homogeneous spray manifold is endowed with a reductive decomposition, we can use the spray vector field to describe these properties, and prove that a weakly symmetric spray manifold must be geodesic orbit, which generalizes its analog in homogeneous Riemannian and Finsler geometries.

\vbox{}
\end{abstract}

\maketitle

\section{Introduction}

A connected Riemannian homogeneous manifold is called {\it geodesic orbit} (or {\it g.o.} in short), if any geodesic is the orbit of a one-parameter subgroup of isometries. This notion was introduced by O. Kowalski and L. Vanhecke in 1991 \cite{KV1991}, as a generalization of naturally reductive homogeneity. Since then, there have been many research works on this subject. See \cite{AA2007, AN2009, BN2020, DKN2004, Go1996, GN2018} and references therein. There are many important subclasses in g.o. Riemannian manifolds, which include the {\it weakly symmetric} Riemannian manifolds. The notion of weak symmetry was proposed by A. Selberg \cite{Se1956} and it has interesting connections with spherical spaces, commutative spaces, Galfand pairs, etc. \cite{AV1999,BV1996,Wo2007,Ya2004,Zi1996}. A connected homogeneous Riemannian manifold is {weakly symmetric} if and only if any two points can be exchanged by an isometry, or equivalently, for any geodesic $c(t)$ there exists an element $g\in G$ such that $g\cdot c(t)=c(-t)$ for all $t\in\mathbb{R}$. It is a well known fact that a weakly symmetric Riemannian manifold is g.o. (see \cite{BN2018,BKV1997} for two different proofs).
Meanwhile, the g.o. and weakly symmetric properties were generalized by S. Deng and his coworkers to Finsler geometry
\cite{DH2010,YD2014} (see \cite{Xu2018,Xu2021,XD2017,ZX2022} for some recent progresses).
By similar arguments as in the Riemannian case, one can prove that a weakly symmetric Riemannian manifold is also g.o. \cite{De2012}.

Beyond Riemann-Finsler geometry, there is a more general spray geometry, where the metrics are ignored and only the (properly parametrized) geodesics are concerned \cite{Sh2001}.
Nowadays, this field has caught the attention of many geometers \cite{LM2021,LS2018,Ya2011},
and it has been studied in the homogeneous context as well \cite{HM2021,Xu2022-1,Xu2022,Xu2023}
Homogeneous spray geometry generalizes homogeneous Finsler geometry. In fact, it generalizes homogeneous pseudo-Riemannian geometry and homogeneous pseudo-Finsler geometry as well.

%In, the notion of geodesic orbit Finsler space was defined, and in \cite{XD2017}, the interaction between geodesic orbit property and negative curvature property was explored.
%\subsection{The definition of weakly symmetric and geodesic orbit in spray geometry}

In this paper, we generalize the g.o. and weakly symmetric properties to homogeneous spray geometry.
To avoid iterations, we automatically assume the homogeneous manifolds are connected and Lie groups have only countable connected components.

\begin{definition}\label{def-spray-g.o. space}
We call a homogeneous spray manifold $(G/H,\mathbf{G})$ {\it g.o.} if each geodesic $c(t)$ with $t$ sufficiently close to $0$ is {\it $G$-homogeneous}, i.e., there exists $X\in \mathfrak{g}=\mathrm{Lie}(G)$, $x\in G/H$ and $\epsilon>0$, such that $c(t)=\exp tX\cdot x$ for $t\in(-\epsilon,\epsilon)$.
\end{definition}

%$(1)$. For any geodesic $c(t)\in (G/H,\mathbf{G})$, after reversing $c(t)$, it is still a geodesic.
%
%$(2)$. For any geodesic $\gamma(t)$ and $p\neq q\in \gamma(t)$, there exists $g\in G$, such that $g(p)=q$, $g(q)=p$, and $g(\gamma)$ is a reversible geodesic.
%\end{definition}

\begin{definition}\label{def-spray-weakly symmetric}
We call a homogeneous spray manifold $(G/H,\mathbf{G})$ {\it weakly symmetric} if for any geodesic $c(t)$ with $t$ sufficiently close to $0$, there exists $g\in G$, such that $g\cdot c(t)=c(-t)$
for all $t$ where both $c(t)$ and $c(-t)$ are defined.
\end{definition}

\begin{remark} A homogeneous spray manifold may not be complete (see Section 4.3 in \cite{Xu2022-1} for an explicit example), so it is not reasonable to assume that the geodesics
$c(t)$ in above definitions are defined for all $t\in\mathbb{R}$.
\end{remark}

Next,
we use the spray vector field $\eta$ \cite{Xu2022-1,Xu2022} to characterize these properties.
Now assume that $G/H$ admits a reductive decomposition $\mathfrak{g}=\mathfrak{h}+\mathfrak{m}$. Then, $\eta$ is an $\mathrm{Ad}(H)$-equivariant positively 2-homogeneous smooth map from $\mathfrak{m}\backslash\{0\}$ to $\mathfrak{m}$. Conversely, any map satisfying this description is the spray vector field for some homogeneous spray structure on $G/H$ (see Theorem 5.5 in \cite{Xu2022}). We usually view $\eta$ as a smooth tangent vector field on $\mathfrak{m}\backslash\{0\}$.

\begin{theorem}\label{thm-1}
Let $(G/H,\mathbf{G})$ be a  homogeneous spray manifold with a reductive decomposition
$\mathfrak{g}=\mathfrak{h}+\mathfrak{m}$ and the spray vector field $\eta$. Then
$(G/H,\mathbf{G})$ is g.o. if and only if $\eta$ is tangent to $\mathrm{Ad}(H)$-orbits in $\mathfrak{m}\backslash\{0\}$, i.e., at each $y\in\mathfrak{m}\backslash\{0\}$, $\eta(y)\in T_{y}(\mathrm{Ad}(H)y)=[\mathfrak{h},y]$.
\end{theorem}

Here the orbit $\mathrm{Ad}(H)y$ for $y\in\mathfrak{m}\backslash\{0\}$ is viewed as a submanifold in $\mathfrak{m}\backslash\{0\}$ (see Theorem 2.9.7 in \cite{Va1984}).
Theorem \ref{thm-1} helps us find g.o. spray manifolds. For example, any homogeneous manifold $G/H$ with a reductive decomposition can be endowed with a canonical g.o. spray structure $\mathbf{G}_0$ with vanishing spray vector field (see page 30 in \cite{Xu2022}).
Theorem \ref{thm-1} also reveals the difference between g.o. properties in Finsler geometry and spray geometry. For example, on $S^{2n-1}=U(n)/U(n-1)$, any homogeneous Finsler metric is g.o. \cite{Xu2018}, but non-g.o. homogeneous spray structures can be easily constructed.

\begin{theorem}\label{thm-2}
Let $(G/H,\mathbf{G})$ be a homogeneous spray manifold with a reductive decomposition
$\mathfrak{g}=\mathfrak{h}+\mathfrak{m}$. Then $(G/H,\mathbf{G})$ is weakly symmetric if and only if for each $y\in\mathfrak{m}\backslash\{0\}$, there exists $g\in H$, such that $\mathrm{Ad}(g)y=-y$. In particular, when $(G/H,\mathbf{G})$ is weakly symmetric, its spray vector field $\eta$ is even, i.e., $\eta(y)=\eta(-y)$, $\forall y\in\mathfrak{m}\backslash\{0\}$.
\end{theorem}

The first statement in Theorem \ref{thm-2} implies that the classification for weakly symmetric spray manifolds is an algebraic problem \cite{WC2022}. Here is a simple example. On any connected Lie group $G$, there is a canonical affine bi-invariant spray structure $\mathbf{G}_0$ \cite{Xu2022-1}, which geodesics are orbits of one-parameter subgroups in $G$. This $\mathbf{G}_0$ is a symmetric (so it is also weakly symmetric) homogeneous spray structure when $G$ is viewed as the homogeneous manifold $\widetilde{G}/\widetilde{H}$, where $\widetilde{G}$ is generated by all left and right translations, and the inverse operation $\tau(g)=g^{-1}$. However, $\mathbf{G}_0$ is g.o. but not weakly symmetric when $G$ is viewed as the homogeneous manifold $G/\{e\}$.

Finally, we use Theorem \ref{thm-1} and Theorem \ref{thm-2} to prove

\begin{theorem}\label{main-theorem}
Any weakly symmetric spray manifold $(G/H,\mathbf{G})$ with a reductive decomposition $\mathfrak{g=h+m}$ is g.o..
\end{theorem}

Theorem \ref{main-theorem} generalizes its analogs in homogeneous Riemannian and Finsler geometries. To prove this theorem, we applied a new approach, which is based on the theoretical framework of homogeneous spray geometry \cite{Xu2022-1,Xu2022,Xu2023}. There is a remaining question:

\begin{question}If a weakly symmetric spray manifold does not admit a reductive decomposition, must it still be g.o.?
\end{question}

%\begin{theorem}\label{main-theorem-3}
%Let $(G/H,\mathbf{G})$ be a homogeneous spray manifold with reductive decomposition $\mathfrak{g=h+m}$, and $\eta$ the spray vector field, $\mathbf{G}$ is $G$ -invariant spray structure. If $(G/H,\mathbf{G})$ is a weakly symmetric homogeneous spray space, then $\eta$ is tangent to $\Ad(H)$ -orbit.
%\end{theorem}

%To use theorem \ref{main-theorem-g.o. spray} and theorem \ref{main-theorem-3}

%\begin{theorem}\label{main-theorem-2}
%Let $(G/H,\mathbf{G})$ be a homogeneous spray manifold with reductive decomposition $\mathfrak{g=h+m}$, and $\eta$ the spray vector field, $\mathbf{G}$ is $G$ -invariant spray structure. If $(G/H,\mathbf{G})$ is a weakly symmetric spray space, then $(G/H,\mathbf{G})$ is a spray g.o. space.
%\end{theorem}

This paper is scheduled as follows.
In Section 2, we summarize some preliminary knowledge in spray geometry and introduce the spray vector field for a homogeneous spray manifold. In Section 3, we prove all the theorems.

\section{Preliminaries}
\subsection{Spray structure and geodesic}

The spray structure on a smooth manifold $M$ is a smooth tangent vector field $\mathbf{G}$ on  $TM\backslash0$, which can be locally presented as
\begin{equation*}
\mathbf{G}=y^{i}\partial_{x^{i}}-2\mathbf{G}^{i}\partial_{y^{i}}
\end{equation*}
for any standard local coordinate $(x,y)=(x^i, y^i)$ (i.e., $x = (x^i) \in M $ and $y = y^j\partial _{x^j} \in T_xM$). Here the spray coefficients $\mathbf{G}^{i}=\mathbf{G}^{i}(x,y)$ are required to be  positively 2-homogeneous for its $y$-entry, i.e., $ \mathbf{G}^{i}(x,\lambda y) = \lambda^{2}\mathbf{G}^{i}(x,y), \forall \lambda > 0$. We call $(M,\mathbf{G})$ a {\it spray manifold}.

A smooth curve $c(t)$ on $(M,\mathbf{G})$ with nowhere-vanishing $\dot{c}(t)$ is called a geodesic if its lifting $ (c(t), \dot{c}(t))$ in $TM\setminus 0$ is an integral
curve of $\mathbf{G}$. Locally, a geodesic $c(t) = (c^i(t))$ can be characterized by the following ODE,
\begin{equation*}
\ddot{c}(t)+2\mathbf{G}^{i}(c(t),\dot{c}(t))=0, \forall i.
\end{equation*}

On a Finsler manifold $(M,F)$, the canonical spray structure $\mathbf{G}_F$
is determined by $$\mathbf{G}^i_F=\tfrac14g^{il}([F^2]_{x^ky^l}y^k-[F^2]_{x^l}).$$
The geodesic and many curvatures (including Riemann curvature, Ricci curvature, S-curvature, etc.) of
$(M,F)$ coincide with those of $(M,\mathbf{G}_F)$ respectively.
See \cite{Sh2001} for more details.

\subsection{Homogeneous spray manifold}
A connected spray manifold $(M,\mathbf{G})$ is called homogeneous if there exists a connected Lie transformation group $G$ which acts transitively on $M$ and preserves the set of all geodesics. More precisely, it means that, for any $g\in G$, its tangent map on $TM$ preserves $\mathbf{G}$ \cite{Xu2022}. In this situation, we may identify $M$ as a homogeneous manifold $G/H$.

We assume that $G/H$ has a {\it reductive} decomposition $\mathfrak{g=h+m}$, where $\mathfrak{g}=\mathrm{Lie}(G)$ and $\mathfrak{h}=\mathrm{Lie}(H)$. The reductiveness here
means that the given decomposition is $\mathrm{Ad}(H)$-invariant, or in the Lie algebraic level,
$[\mathfrak{h},\mathfrak{m}] \subset\mathfrak{m}$. Then the subspace $\mathfrak{m}$ can be equivariantly identified with the tangent space $T_o(G/H)$ at $o=eH$. Here
the equivariance means that
the $\mathrm{Ad}(H)$-action on $\mathfrak{m}$ coincides with
the isotropy action on $T_o(G/H)$.

\begin{remark}
A homogeneous Finsler manifold always has a reductive decomposition \cite{De2012,Ni2017}, but a homogeneous spray manifold may not.
\end{remark}
\subsection{Spray vector field}
Let $(G/H,\mathbf{G})$ be a homogeneous spray manifold with a reductive decomposition
$\mathfrak{g}=\mathfrak{h}+\mathfrak{m}$.
The {\it spray vector field} $\eta:\mathfrak{m}\backslash\{0\}\rightarrow\mathfrak{m}$ for $(G/H,\mathbf{G})$ is defined as follows \cite{Xu2022}. Let $c(t)$
be a geodesic on $(G/H,\mathbf{G})$ satisfying $c(0)=o$ and $\dot{c}(0)=y\in T_o(G/H)\backslash\{0\}=\mathfrak{m}\backslash\{0\}$.
By the following lemma (Lemma 5.1 in \cite{Xu2022}),
we can uniquely lift $c(t)$ to a curve in $G$.

\begin{lemma}\label{lemma5.1}
There exists a unique smooth curve $\overline{c}(t)$ on $G$ satisfying $\overline{c}(0)=e$, $c(t)=\overline{c}(t)\cdot o$,  and  $(L_{\overline{c}(t)^{-1}})_{*}(\dot{\overline{c}}(t))=
(\overline{c}(t)^{-1})_{*}(\dot{c}(t))\in\mathfrak{m}\setminus\{0\}$ when $t$ is sufficiently close to $0$.
\end{lemma}

\begin{remark} For simplicty, we will always assume the parameter $t$ in later discussions to be sufficiently close to $0$ without repeating.
\end{remark}
%Further more, $\overline{c}(t)$ can also be defined for $t\in (a,b)$ if one of the following is satisfied:
%\par
%$(1)$. The decomposition $\mathfrak{g = h + m}$ is reductive;
%\par
%$(2)$. The subgroup $H$ is compact.
%\end{lemma}
%The proof of this lemma is important for another theorem in this paper. In proof of this lemma, we get a useful equation.
%\begin{align}
%(L_{\overline{c}(t)^{-1}})_{*}(\dot{\overline{c}}(t))=(L_{h(t)^{-1}})_{*}(\dot{h}(t))+Ad(h(t)^{-1})(L_{g(t)^{-1}})_{*}(\dot{g}(t)).
%\end{align}
%The specific details of the proof can be found in \cite{Xu2022}.

Denote by $y(t)=(L_{\overline{c}_{y}(t)^{-1}})_{*}(\dot{\overline{c}}_{y}(t))$ a smooth curve in  $\mathfrak{m}\setminus\{0\}$, then we define the spray vector field $\eta$ by $\eta(y)= -\frac{d}{dt}|_{t=0}y(t)$.
Theorem 5.3 in \cite{Xu2022} is crucial for proving our theorems later. We reformulate it as follows.

\begin{theorem}\label{def-eta}
Let $(G/H,\mathbf{G})$ be a connected homogeneous spray manifold with a reductive decomposition $\mathfrak{g=h+m}$ and the spray vector field $\eta$. Then there is a one-to-one correspondence between the following two sets:
\begin{enumerate}
\item The set of all geodesics $c(t)$ with $c(0) = o$;
\item The set of all integral curves $y(t)$ of $-\eta$ which is defined at $t=0$.
\end{enumerate}
The correspondence is from $c(t)$ to $y(t) = (L_{\overline{c}(t)^{-1}})_{*}(\dot{\overline{c}}(t))$,
where $\overline{c}(t)$ is the lifting of $c(t)$ provided by Lemma \ref{lemma5.1}.
\end{theorem}

We will also needs the following property of $\eta$ (see Lemma 5.2 and Theorem 5.5 in \cite{Xu2022}).

\begin{lemma}\label{lemma-2}
Let $(G/H,\mathbf{G})$ be a connected homogeneous spray manifold with a reductive decomposition $\mathfrak{g=h+m}$, then its spray vector field is $\mathrm{Ad}(H)$-equivariant, i.e.,
$\forall g\in H$, $y\in\mathfrak{m}\backslash\{0\}$, $\eta(\mathrm{Ad}(g)y)=\mathrm{Ad}(g)\eta(y)$.
\end{lemma}

\begin{remark}Spray vector field was first introduced by L. Huang for a homogeneous Finsler manifold
with a reductive decomposition \cite{Hu2015}.
When we generalize this notion to homogeneous spray geometry, reductive decomposition is not essential. All above statements except the $\mathrm{Ad}(H)$-equivariance of $\eta$ are still valid.
\end{remark}
%\section{The properties of symmetric and geodesic orbit on spray structure}%Иэёц¶ЁАнµДРрКцј°Ц¤Гч
%In subsection Finsler g.o.space, we introduce the definition \ref{def-g.o.space} of Finsler g.o.space and definition \ref{def-weakly symmetric}, now we can give a equivalent description of g.o. space and weakly symmetric space with spray structure. Firstly, let we introduce the equivalence theorem of g.o. space with spray structure.

\section{Proofs of the main Results}

\subsection{Proof of Theorem \ref{thm-1}}
Let $(G/H,\mathbf{G})$ be a homogeneous spray manifold with a reductive decomposition
$\mathfrak{g=h+m}$ and the spray vector field $\eta$.
First, we assume the g.o. property and prove that $\eta$ is tangent to the $\mathrm{Ad}(H)$-orbits in $\mathfrak{m}\backslash\{0\}$.

We choose a vector $y$ from $\mathfrak{m}\backslash\{0\}$ and denote by
$c(t)$ be the geodesic on $(G/H,\mathbf{G})$ which satisfies $c(0)=o$ and $\dot{c}(0)=y\in\mathfrak{m}\backslash\{0\}$. Then the g.o. property provides a vector $u\in\mathfrak{g}$, such that
$c(t)=\exp tu\cdot o$. Notice that $u-y\in\mathfrak{h}$.
Lemma \ref{lemma5.1} provides a smooth curve $\overline{c}(t)$ in $G$ such that
\begin{eqnarray}
& &\overline{c}(0)=e,\quad c(t)=\exp tu\cdot o=\overline{c}(t)\cdot o,\label{001}\\
& &(L_{\overline{c}(t)^{-1}})_*(\dot{\overline{c}} (t))\in\mathfrak{m}\backslash\{0\}. \label{002}
\end{eqnarray}
By (\ref{001}), $h(t)=\exp (-tu)\overline{c}(t)$
is a smooth curve in $H$ satisfying $h(0)=e$.
Because $$\dot{\overline{c}}(t)=\tfrac{d}{dt}(\exp (tu) h(t))=(L_{\exp tu})_*\dot{h}(t)+(R_{h(t)})_*(L_{\exp tu})_*u,$$
we have
\begin{eqnarray}
(L_{\overline{c}(t)^{-1}})_*(\dot{\overline{c}} (t))&=&
(L_{h(t)^{-1}})_*[(L_{\exp(-tu)})_*(L_{\exp tu})_*\dot{h}(t)+(R_{h(t)})_*u]\nonumber\\
&=&(L_{h(t)^{-1}})_*\dot{h}(t)+\mathrm{Ad}({h(t)^{-1}})u\nonumber\\
&=&[(L_{h(t)^{-1}})_*\dot{h}(t)+\mathrm{Ad}({h(t)^{-1}})(u-y)]+
[\mathrm{Ad}({h(t)^{-1}})y].\label{003}
\end{eqnarray}
By the reductiveness of the given decomposition, the two summands in right side of  (\ref{003})
belong to $\mathfrak{h}$ and $\mathfrak{m}$ respectively. So (\ref{002}) implies
$y(t)=(L_{\overline{c}(t)^{-1}})_*(\dot{\overline{c}} (t))=\mathrm{Ad}({h(t)^{-1}})y$. So $y(t)$ is a smooth curve on $\mathrm{Ad}(H)y$, and then its tangent vector
$\tfrac{d}{dt}|_{t=0}y(t)=[y,\dot{h}(0)]=-\eta(y)$ is contained in $T_{y}(\mathrm{Ad}(H)y)=[\mathfrak{h},y]$.
One side of Theorem \ref{thm-1} is proved.

Next, we assume that $\eta$ is tangent to $\mathrm{Ad}(H)$-orbits in $\mathfrak{m}\backslash\{0\}$ and prove the g.o. property. Choose any vector $y$ from $\mathfrak{m}\backslash\{0\}$, then our assumption provides $\eta(y)=[v,y]$ for some $v\in\mathfrak{h}$.
%Denote by $y(t)$ the integral curve
%of $-\eta$ which satisfies $y(0)=y$.

\noindent{\bf Claim A}:
$y(t)=\mathrm{Ad}(\exp(-tv))y$ is the integral curve of $-\eta$ which satisfies $y(0)=y$.

It is obvious that $y(0)=y$.
By Lemma \ref{lemma-2},
\begin{equation}\label{006}
\eta(y(t))=\mathrm{Ad}(\exp(-tv))\eta(y)=
\mathrm{Ad}(\exp(-tv))[v,y]=[u,\mathrm{Ad}(\exp(-tv))y].
\end{equation}
On the other hand,
\begin{equation}\label{007}
\dot{y}(t)=\tfrac{d}{dt}e^{-t\mathrm{ad}(v)}y=-\mathrm{ad}(v)e^{-t\mathrm{ad}(v)}y=
-[v,\mathrm{Ad}(\exp(-tv))y].
\end{equation}
Comparing (\ref{006}) and (\ref{007}), we get $\dot{y}(t)=-\eta(y(t))$. Claim A is proved.

Let $c(t)$ be the geodesic on $(G/H,\mathbf{G})$ which satisfies $c(0)=o$ and $\dot{c}(0)=y$. Then Theorem \ref{def-eta} tells us $y(t)=(L_{\overline{c}(t)^{-1}})_*(\dot{\overline{c}}(t))$.
Denote by $g(t)$ the curve in $G$ determined by $\overline{c}(t)=g(t)\exp (tv)$. Obviously $g(0)=0$. Similar
calculation as in (\ref{003}) or in the proof of Lemma 5.1 in \cite{Xu2022} shows
\begin{eqnarray}
(L_{\overline{c}(t)^{-1}})_*(\dot{\overline{c}}(t))
&=&
(L_{\exp(-tv)})_*\tfrac{d}{dt}\exp tv+\mathrm{Ad}(\exp(-tv))((L_{g(t)^{-1}})_*\dot{g}(t))\nonumber\\
&=&v+\mathrm{Ad}(\exp(-tv))((L_{g(t)^{-1}})_*\dot{g}(t))\nonumber\\
&=&[v+\mathrm{Ad}(\exp(-tv))(((L_{g(t)^{-1}})_*\dot{g}(t))_\mathfrak{h})]\nonumber\\
& &+
[\mathrm{Ad}(\exp(-tv))(((L_{g(t)^{-1}})_*\dot{g}(t))_\mathfrak{m})],\label{004}
\end{eqnarray}
in which the subscripts $\mathfrak{h}$ and $\mathfrak{m}$ mean the projection according to the given reductive decomposition.
The two summands in the right side of (\ref{004})
belongs to $\mathfrak{h}$ and $\mathfrak{m}$ respectively.
Combine Claim A and (\ref{004}), we see
\begin{equation*}
v+\mathrm{Ad}(\exp(-tv))(((L_{g(t)^{-1}})_*\dot{g}(t))_\mathfrak{h})=0\ \mbox{and}\
((L_{g(t)^{-1}})_*\dot{g}(t))_\mathfrak{m}=y,
\end{equation*}
i.e., $(L_{g(t)^{-1}})_*\dot{g}(t)\equiv y-v$.
So $g(t)=\exp t (y-v)$ and $c(t)=\overline{c}(t)\cdot = g(t)\exp tv \cdot o=\exp t(y-v)\cdot o$.
We get the g.o. property for $(G/H,\mathbf{G})$ and
end the proof of Theorem \ref{thm-1}.

\subsection{Proof of Theorem \ref{thm-2}}
Let $(G/H,\mathbf{G})$ be a weakly symmetric spray manifold with a reductive decomposition
$\mathfrak{g=h+m}$.
We choose a vector $y$ from $\mathfrak{m}\backslash\{0\}$, and denote by $c_1(t)$ the geodesic on $(G/H,\mathbf{G})$ which satisfies $c(0)=o$ and $\dot{c}(0)=y$.
By Definition \ref{def-spray-weakly symmetric}, there exists $g\in H$, such that the geodesic $c_2(t)=g\cdot c(t)$ coincides with $c_1(-t)$. Obviously $\dot{c}_2(0)=\mathrm{Ad}(g)y=-\dot{c}_1(0)=-y$. This argument can be reversed, which proves
the first statement in Theorem \ref{thm-2}.

Let $\overline{c}_{1}(t)$ be the lifting of $c_{1}(t)$
provided by Lemma \ref{lemma5.1}.

\noindent{\bf Claim B}:
$\overline{c}_{2}(t)=\overline{c}_{1}(-t)$ is the lifting of $c_2(t)$ provided by Lemma \ref{lemma5.1}.

It is obvious that $\overline{c}_{2}(0)=\overline{c}_{1}(0)=e$ and
$
c_2(t)=c_1(-t)=\overline{c}_1(-t)\cdot o=\overline{c}_2(t)\cdot o$.  Meanwhile, we have
$\dot{\overline{c}}_2(t)=\tfrac{d}{dt}\overline{c}_1(-t)=-\dot{\overline{c}}_1(-t)$, and then
\begin{eqnarray*}
(L_{\overline{c}_2(t)^{-1}})_*(\dot{\overline{c}}_2(t))=(L_{\overline{c}_1(-t)^{-1}})_*(-\dot{\overline{c}}_1(-t))
=-(L_{\overline{c}_1(-t)^{-1}})_*(\dot{\overline{c}}_1(-t))\in\mathfrak{m}\backslash\{0\}.
\end{eqnarray*}
So $\overline{c}_2(t)=\overline{c}_1(-t)$ is the lifting of $c_2(t)=c_1(-t)$ provided by Lemma \ref{lemma5.1}, which proves Claim B.

By the definition of $\eta$, we have $$\eta(y)=-\tfrac{d}{dt}|_{t=0}[(L_{\overline{c}_{1}(t)^{-1}})_{*}\dot{\overline{c}}_{1}(t)]\quad
\mbox{and}\quad  \eta(-y)=-\tfrac{d}{dt}|_{t=0}[(L_{\overline{c}_{1}(t)^{-1}})_{*}\dot{\overline{c}}_{1}(t)].$$ Direct calculation shows
\begin{eqnarray*}
\eta(-y)&=&
-\tfrac{d}{dt}|_{t=0}[(L_{\bar{c}_{2}(t)^{-1}})_{*}(\dot{\bar{c}}_{2}(t))]
%&=-\tfrac{d}{dt}|_{t=0}[(L_{\bar{c}_{1}(-t)^{-1}})_{*}(\tfrac{d}{dt}{\bar{c}}_{1}(-t))] \\
=
\tfrac{d}{dt}|_{t=0}[(L_{\bar{c}_{1}(-t)^{-1}})_{*}(\dot{\bar{c}}_{1}(-t)] \\
%&=-\frac{d}{d(-t)}[(L_{\bar{c}_{1}(-t)^{-1}})_{*}\dot{\bar{c}}_{1}(-t)]|_{-t=0}\notag\\
&=&-\tfrac{d}{ds}|_{s=-t=0}[(L_{\bar{c}_{1}(s)^{-1}})_{*}\dot{\bar{c}}_{1}(s)]
=\eta(y),
\end{eqnarray*}
which proves Theorem \ref{thm-2}.

\subsection{Proof of Theorem \ref{main-theorem}}
Let $(G/H,\mathbf{G})$ be a weakly symmetric spray manifold with a reductive decomposition $\mathfrak{g=h+m}$ and the spray vector field $\eta$.
We choose any vector $y$ from $\mathfrak{m}\backslash\{0\}$ and denote by $y_1(t)$ the integral curve of $-\eta$ which satisfies $y_1(0)=y$. By Theorem \ref{thm-1}, to prove Theorem \ref{main-theorem}, we only need to prove that, for each value $t_1$, $y_1(t_1)$ is contained in $\mathrm{Ad}(H)y$. It should be notified that, we have used here the quasi regularity of the submanifold $\mathrm{Ad}(H)y$ in $\mathfrak{m}\backslash\{0\}$ (see page 117 in \cite{GHV1973}).

By the first statement in Theorem \ref{thm-2}, there exists $g_1,g_2\in H$ satisfying $\mathrm{Ad}(g_1)y=-y$ and $\mathrm{Ad}(g_2)y_1(\tfrac{t_1}2)=-y_1(\tfrac{t_1}2)$. Denote by $y_2(t)=
\mathrm{Ad}(g_2) y_1(t+\tfrac{t_1}2)$ another curve in $\mathfrak{m}\backslash\{0\}$.
Calculation shows
\begin{eqnarray*}
\dot{y}_2(t)&=&\tfrac{d}{dt}(\mathrm{Ad}(g_2)y_1(t+\tfrac{t_1}2))
=\mathrm{Ad}(g_2)\dot{y}_1(t+\tfrac{t_1}2)
=\mathrm{Ad}(g_2)(-\eta(y_1(t+\tfrac{t_1}2)))\\
&=&-\eta(\mathrm{Ad}(g_2)y_1(t+\tfrac{t_1}2))=-\eta(y_2(t)),
\end{eqnarray*}
in which Lemma \ref{lemma-2} has been used for the fourth equality.
It indicates that $y_2(t)$ is the integral curve of $-\eta$ which passes $-y_1(\tfrac{t_1}2)$ at $t=0$.
On the other hand,
we have
\begin{equation*}\label{010}
\tfrac{d}{dt}(-y_1(-t+\tfrac{t_1}2))
=\tfrac{d}{ds}|_{s=-t+\tfrac{t_1}2}y_1(s)=-\eta(y_1(-t+\tfrac{t_1}2))=-\eta(-y_1(t+\tfrac{t_1}2)),
\end{equation*}
where the second statement in Theorem \ref{thm-2} has been used for the last equality.
It tells us that $-y_1(t+\tfrac{t_1}2)$ is also the integral curve of $-\eta$ which passes $-y_1(\tfrac{t_1}2)$ at $t=0$.
So we must have
\begin{equation}\label{011}
\mathrm{Ad}(g_2)y_1(t+\tfrac{t_1}2)\equiv-y_1(-t+\tfrac{t_1}2).
\end{equation} 
Input $t=\tfrac{t_1}2$ into (\ref{011}),  then we get
$\mathrm{Ad}(g_2)y_1(t_1)=-y_1(0)=-y=\mathrm{Ad}(g_1)y$, i.e., $y_1(t_1)\in\mathrm{Ad}(H)y$. The proof of Theorem \ref{main-theorem} is finished.
\bigskip

\noindent
{\bf Acknowledgement}\quad
This paper is supported by National Natural Science Foundation of China (No. 11821101, No. 12131012), and Beijing Natural Science Foundation (No. 1222003).

\end{document}